\documentclass[psamsfonts]{amsart}

\usepackage{amssymb,amsfonts}
\usepackage[all,arc]{xy}
\usepackage{enumerate}
\usepackage{mathrsfs}
\usepackage{setspace}
\usepackage{graphicx}
\usepackage{caption}
\usepackage{float}
\DeclareGraphicsExtensions{.pdf,.png,.jpg,.eps}

\theoremstyle{definition}

\theoremstyle{remark}

\makeatletter
\let\c@equation\c@thm
\makeatother
\numberwithin{equation}{section}
\graphicspath{ {desktop/} }

\title{An Extension of Heron's Formula}

\author{Zohreh Shahbazi}

\date{March 15, 2017}

\begin{document}

\begin{abstract}

This paper introduces an extension of Heron's formula to approximate area of cyclic n-gons where the error never exceeds $\frac{\pi}{e}-1$.
\end{abstract}
\maketitle
\section{Introduction}

A cyclic n-gon is a polygon with $n$ vertices all located on the perimeter of a given circle. Triangles and cyclic quadrilaterals are examples of cyclic n-gons. We are interested in finding a formula for the area of a cyclic $n$-gon in terms of its sides. For $n=3$ and $n=4$ the formula is known, namely Heron's and Brahmagupta's formulas.
Heron's formula calculates the area of a given triangle in terms of its three sides: $$S=\sqrt{P(P-a)(P-b)(P-c)}$$ where $a, b, c$ are sides, $P=\frac{1}{2}(a+b+c)$ is the semi-perimeter, and $S$ is the area of the triangle. Brahmagupta's formula calculates the area of a cyclic quadrilateral in terms of its four sides: $$S=\sqrt{(P-a)(P-b)(P-c)(P-d)}$$ where $a, b, c, d$ are sides, $P=\frac{1}{2}(a+b+c+d)$ and $S$ is the area of the cyclic quadrilateral.\par

Naturally we can ask if there is any similar formula to calculate area of a cyclic $n$-gon in terms of its sides. In this paper, we will introduce a natural extension of Heron and Brahmagupta's formulas. It turns out that this generalized formula is not exactly equal to the area of the cyclic n-gon, but it can approximate the area with a small error which never exceeds $\frac{\pi}{e}-1$.

The area of a general convex quadrilateral could be obtained by formulas which have terms similar to Heron's terms \cite{Jos}. It is also known that the area of a cyclic n-gon raised to power two and then multiplied  with a factor sixteen, is a monic polynomial whose other coefficients are polynomials in the sides of n-gon \cite{igor}, \cite{rob}, \cite{v}. Therefore, the formula suggested in this paper of area an n-gon in terms of its sides provides with an estimation for Heron polynomials \cite{connelly}.

\section{Area of n-gons in terms of its sides}
Imagine a cyclic n-gon with sides $x_1, x_2, ...,x_n$, area $S_n$ and semi-perimeter $$P_n=\frac{1}{2}(x_1+x_2+...+x_n).$$ A natural extension of Heron or Brahmagupta's formula is $$S_n=\sqrt{P^{(4-n)}(P-x_1)(P-x_2)...(P-x_n)}.$$ Here $n\geq 3$. Obviously, when $n=3$, we get Heron's formula and when $n=4$, we get Brahmagupta's. Another good aspect of this suggested formula is that if the length of one of the sides of n-gon, say $x_1=0$, then we obtain an $(n-1)$-gon. The new formula consistently calculates the area of resulting $(n-1)$-gon in two ways:$$S_n=\sqrt{P^{(4-n)}(P-0)(P-x_2)...(P-x_n)}=\sqrt{P^{(4-(n-1))}(P-x_2)...(P-x_n)}=S_{n-1}.$$

To understand if this formula is correct to calculate the area of n-gon, we first compute the area in terms of central angles: \\
Connect all vertices to the center of circle, which inscribes the $n$-gon. Denote the central angles $\alpha_1, \alpha_2, ..., \alpha_n$. So, $\sum_{i=1}^n \alpha_i=2\pi$ (Figure 1.)
\begin{figure}[H]
	\begin{center}
		\includegraphics[width = 0.5\textwidth, height = 0.5\textwidth]{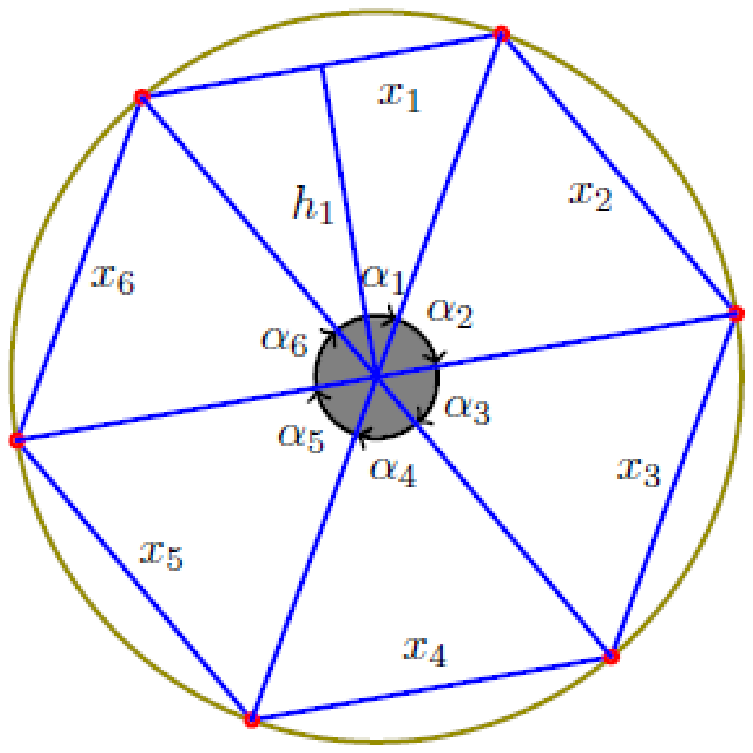}
	\end{center}	
\caption{}
\end{figure}
Denote the exact area of the cyclic n-gon by $A_n$. So,
\begin{equation*}
\begin{aligned}
A_n &=\frac{1}{2}\sum_{i=1}^n h_ix_i\\
&= \frac{1}{4}\sum_{i=1}^n \cot \frac{\alpha_i}{2}x_i^2\\
&=\frac{1}{4}\sum_{i=1}^n 4 R^2 \sin^2 \frac{\alpha_i}{2}\cot \frac{\alpha_i}{2}\\
&=R^2\sum_{i=1}^n\sin \frac{\alpha_i}{2}\cos \frac{\alpha_i}{2}\\
&=\frac{R^2}{2}\sum_{i=1}^n\sin \alpha_i.\\
\end{aligned}
\end{equation*}
Here, we use the fact that $h_i= \frac{1}{2}x_i\cot \frac{\alpha_i}{2}$ and $x_i=2R \sin \frac{\alpha_i}{2}$ where $R$ is the radius of the circle.
Next, we calculate generalized Heron's formula for the area of a cyclic n-gon in terms of central angles:
\begin{equation*}
\begin{aligned}
S_n &=\sqrt{P^{(4-n)}(P-x_1)(P-x_2)...(P-x_n)}\\
&= P^2\sqrt{(1-\frac{x_1}{P})(1-\frac{x_2}{P})...(1-\frac{x_n}{P})}\\
&=P^2\prod_{i=1}^n [(1-\frac{x_i}{P})]^{\frac{1}{2}}\\
&=\frac{1}{4}(\sum_{i=1}^nx_i)^2\prod_{i=1}^n [(1-\frac{2x_i}{(\sum_{i=1}^nx_i)})]^{\frac{1}{2}}\\
&=\frac{1}{4}.4.R^2(\sum_{i=1}^n\sin \frac{\alpha_i}{2})^2\prod_{i=1}^n [(1-\frac{2\sin \frac{\alpha_i}{2}}{(\sum_{i=1}^n\sin \frac{\alpha_i}{2})})]^{\frac{1}{2}}.\\
&=R^2(\sum_{i=1}^n\sin \frac{\alpha_i}{2})^2\prod_{i=1}^n [(1-\frac{2\sin \frac{\alpha_i}{2}}{(\sum_{i=1}^n\sin \frac{\alpha_i}{2})})]^{\frac{1}{2}}.\\
\end{aligned}
\end{equation*}
Now we define a new function $D_n$ by dividing the exact area formula $A_n$ from the suggested area formula $S_n$:
\begin{equation*}
\begin{aligned}
D_n &= \frac{S_n}{A_n}\\
&=\frac{R^2(\sum_{i=1}^n\sin \frac{\alpha_i}{2})^2\prod_{i=1}^n [(1-\frac{2\sin \frac{\alpha_i}{2}}{(\sum_{i=1}^n\sin \frac{\alpha_i}{2})})]^{\frac{1}{2}}}{(\frac{R^2}{2}\sum_{i=1}^n\sin \alpha_i)}\\
&= \frac{2(\sum_{i=1}^n\sin \frac{\alpha_i}{2})^2\prod_{i=1}^n [(1-\frac{2\sin \frac{\alpha_i}{2}}{(\sum_{i=1}^n\sin \frac{\alpha_i}{2})})]^{\frac{1}{2}}}{(\sum_{i=1}^n\sin \alpha_i)}.\\
\end{aligned}
\end{equation*}
When all $\alpha_i$, except four, are zero the n-gon is in fact a cyclic quadrilateral and thus $D_n=1$ by Brahmagupta's formula. Let's try to find critical values of $D_n$ as a function of n- variables $\alpha_1, \alpha_2,...,\alpha_n$  with respect to the condition  $\sum_{i=1}^n \alpha_i=2\pi$. By the
Lagrange Multiplier method, this function will be optimized when $$\frac{\partial D_n}{\partial \alpha_1}=\frac{\partial D_n}{\partial \alpha_2}=...=\frac{\partial D_n}{\partial \alpha_n}.$$ Since $D_n$ is a symmetric function with respect to all $\alpha_i$, the solution of system of differential equations will happen at $\alpha_1= \alpha_2=...=\alpha_n$.  

Denote $D_n$ with $D^{reg}_n$ when  all $\alpha_i$ are equal. This will give us the value of $D_n$ for a regular $n$-gon:\\
\begin{equation*}
\begin{aligned}
D^{reg}_n &= \frac{S^{reg}_n}{A^{reg}_n}\\
&=\frac{\frac{n^2x^2}{4} (1-\frac{2}{n})^\frac{n}{2}}{\frac{nx^2}{4}\cot\frac{\pi}{n}}\\
&= n(\tan\frac{\pi}{n})(1-\frac{2}{n})^\frac{n}{2}.\\
\end{aligned}
\end{equation*}
where here $x$ is the length of each side of the regular $n$-gon. Consider the sequence $\{x_n\}_{n=3}^{n=\infty}$ where $ x_n=n(\tan\frac{\pi}{n})(1-\frac{2}{n})^\frac{n}{2}$. The first two terms $x_3$ and $x_4$ are equal 1 by Heron's and Brahmagupta's formulas but $x_5=1.013$. Therefore, equal values of $\alpha_i$ must maximize $D_5$. This proves that the sequence is increasing for $n>4$.  The reason for this behavior is that an $(n-1)$-gon could be considered as a limit case of an $n$- gon, and therefore the maximum value of $D_n$ must be bigger than of maximum value of $D_{n-1}$. So, equal values of $\alpha_i$ must maximize $D_n$. Also, $$\lim_{n\rightarrow \infty} x_n= \lim_{n\rightarrow \infty}(n\tan\frac{\pi}{n})\lim_{n\rightarrow \infty}(1-\frac{2}{n})^\frac{n}{2}=\pi. \frac{1}{e}=\frac{\pi}{e}.$$ This shows that our suggested extension of Heron's formula approximates area of cyclic n-gons where the error never exceeds $\frac{\pi}{e}-1\approx 0.1557.$ In Figure 2, we can see the curve of the function $$f(x)=x\tan\frac{\pi}{x}(1-\frac{2}{x})^\frac{x}{2}$$ for $x\geq 2$. 
\begin{figure}[H]
\begin{center}
\includegraphics[width = 0.99\textwidth, height = 0.3\textwidth]{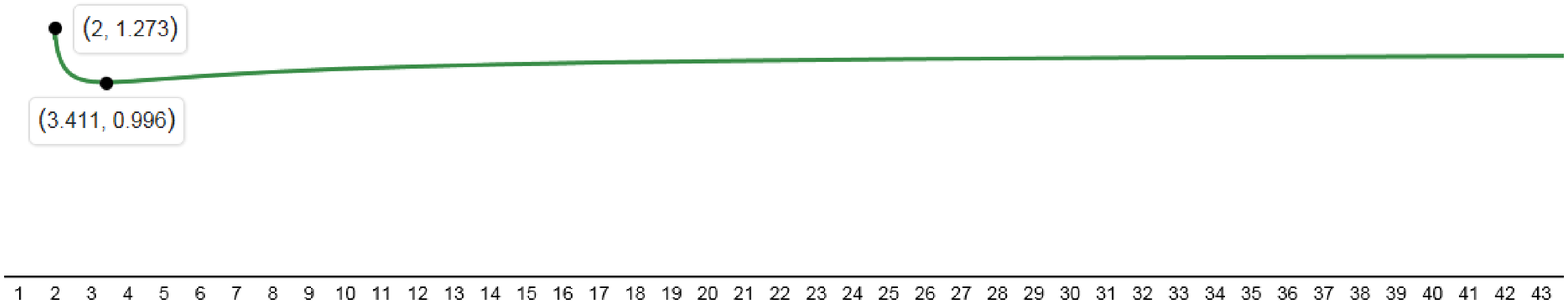}
\end{center}
\caption{}
\end{figure}

\section{Comments}
The area of a general convex quadrilateral could be obtained by formulas which have terms similar to Heron's terms \cite{Jos}. It is also known that the area of a cyclic n-gon raised to power two and then multiplied  with a factor sixteen, is a monic polynomial whose other coefficients are polynomials in the sides of n-gon \cite{igor}, \cite{rob}, \cite{v}. Therefore, the formula suggested in this paper of area an n-gon in terms of its sides provides with an estimation for Heron polynomials \cite{connelly}.


\begin{thebibliography}{1}
	
	\bibitem{Jos}

M.  Josefsson, \emph{Heron-like formulas for quadrilaterals}, The Mathematical Gazette Volume 100, Issue 549 (2016),

505--508. 
\bibitem{igor}
	
	Igor Pak, \emph{The area of cyclic polygons: recent progress on Robbins' Conjecture}, Adv. in Appl. Math 34(4) (2005),
	
	690--696. MR2006b:51017	 
	

	\bibitem{rob}
	
	D. P. Robbins, \emph{Areas of polygons inscribed in a circle}, Discrete Comput. Geom.   12(2) (1994),
	
	223--236. MR95g:51027
	\bibitem{v}
	
	V. V. Varfolomeev, \emph{Inscribed polygons and Heron polynomials (Russian. Russian summary)}, Mat. Sb. 194(3) (2003),
	
	3--24. MR2004d:51014
	
	\bibitem{connelly}

Robert Connelly Pak, \emph{Comments on generalized Heron Polynomials and Robbins' conjectures}, Discrete Mathematics 309 (2009), 4192-4196	
\end{thebibliography}
\end{document}